\documentclass{article}
\usepackage[english]{babel}
\usepackage{latexsym,amssymb,amsmath}
\usepackage{amsfonts,amssymb}
\usepackage{euscript}
\newcounter{theorem} 
\newcounter{lemma} 
\renewcommand{\thetheorem}{\arabic{theorem}}
\renewcommand{\thelemma}{\arabic{lemma}}
\newcommand{\theor}{\par\refstepcounter{theorem}%
{\bf Theorem \thetheorem .}\,\,}

\newcommand{\lem}{\par\refstepcounter{lemma}%
{\bf Lemma \thelemma .}\,\,}
\def\Re{\mathop\mathrm{Re}}\,
\def\Im{\mathop\mathrm{Im}}\,

\begin{document}

\noindent {\bf \Large Schwartz-type integrals in a biharmonic
plane}
\vspace{2mm}\\
\noindent {\bf \large S.V.~Gryshchuk, S.A.~Plaksa}\\
\vspace{5mm}

\begin{abstract}
We consider the commutative algebra $\mathbb{B}$ over the field of
complex numbers with the bases $\{e_1,e_2\}$ satisfying the
conditions $(e_1^2+e_2^2)^2=0$, $e_1^2+e_2^2\ne 0$. The algebra
$\mathbb{B}$ is associated with the biharmonic equation. For
monogenic functions with values in $\mathbb{B}$, we consider a
Schwartz-type boundary value problem (associated with the main
biharmonic problem) for a half-plane and for a disk of the
biharmonic plane $\{xe_1+ye_2\}$, where $x,y$ are real. We obtain
solutions in explicit forms by means of Schwartz-type integrals
and prove that the mentioned problem is solvable unconditionally
for a half-plane but it is solvable for a disk if and only if a
certain natural condition is satisfied.
\end{abstract}

 \vspace{5mm}

\section{Introduction}

We say that an associative commutative two-dimensional algebra
$\mathbb B$ with the unit $1$ over the field of complex numbers
$\mathbb C$ is {\it biharmonic} if in $\mathbb B$ there exists a
{\it biharmonic} basis, i.e a bases $\{e_1,e_2\}$ satisfying the
conditions
\begin{equation}\label{biharm-bas}
 (e_1^2+e_2^2)^2=0,\qquad e_1^2+e_2^2\ne 0\,.
\end{equation}

V.~F. Kovalev and I.~P. Mel'nichenko \cite{Kov-Mel}
 found a
multiplication table for a biharmonic basis $\{e_1,e_2\}$:
\begin{equation} \label{tab_umn_bb}
e_1=1,\qquad 
e_2^2=e_1+2ie_2, 
\end{equation}
where $i$\,\, is the imaginary complex unit. In the paper
\cite{Mel86} I.~P. Mel'nichenko proved that there exists the
unique biharmonic algebra $\mathbb{B}$ with a non-biharmonic basis
$\{1,\rho\}$ for which
\begin{equation} \label{rho}
\rho=2e_1+2ie_2
\end{equation}
and $\rho^{2}=0$, and he constructed all biharmonic bases in
$\mathbb{B}$.

Consider a {\it biharmonic plane} $\mu:=\{\zeta=x\,e_1+y\,e_2 :
x,y\in\mathbb R\}$ which is a linear span of the elements
$e_1,e_2$ of the biharmonic basis (\ref{tab_umn_bb}) over the
field of real numbers $\mathbb R$. With a domain $D$ of the
Cartesian plane $xOy$ we associate the congruent domain
$D_{\zeta}:= \{\zeta=xe_1+ye_2 : (x,y)\in D\}$ in the biharmonic
plane $\mu$. In what follows, $\zeta=x\,e_1+y\,e_2$ and
$x,y\in\mathbb{R}$.

Inasmuch as divisors of zero don't belong to the biharmonic plane,
one can define the derivative $\Phi'(\zeta)$ of function $\Phi :
D_{\zeta}\longrightarrow \mathbb{B}$ in the same way as in the
complex plane:
\[\Phi'(\zeta):=\lim\limits_{h\to 0,\, h\in\mu}
\bigl(\Phi(\zeta+h)-\Phi(\zeta)\bigr)\,h^{-1}\,.\]
 We say that a function $\Phi : D_{\zeta}\longrightarrow \mathbb{B}$ is
\textit{monogenic} in a domain $D_{\zeta}$ if the derivative
$\Phi'(\zeta)$ exists in every point $\zeta\in D_{\zeta}$.

It is established in the paper \cite{Kov-Mel} that a function $\Phi :
D_{\zeta}\longrightarrow \mathbb{B}$ is monogenic in a domain $D_{\zeta}$ if
and only if the following Cauchy--Riemann condition is satisfied
\begin{equation}\label{usl_K_R}
\frac{\partial \Phi(\zeta)}{\partial y}=\frac{\partial \Phi(\zeta)}{\partial x}\,e_2.
\end{equation}

It is proved in the paper \cite{Kov-Mel} that a function
$\Phi(\zeta)$ having derivatives till fourth order in $D_{\zeta}$
satisfies the two-dimensional biharmonic equation
\begin{equation}\label{big-eq}
(\Delta_2)^{2}U(x,y):= \left(\frac{\partial^4}{\partial
x^4}+2\,\frac{\partial^4}{\partial x^2\partial
y^2}+\frac{\partial^4}{\partial y^4}\right)U(x,y)=0
\end{equation}
in the domain $D$ owing to the relations (\ref{biharm-bas}) and
\[(\Delta_2)^{2}\Phi(\zeta)=\Phi^{(4)}(\zeta)\,(e_1^2+e_2^2)^2.\]
Therefore, every component $U_{k}\colon D\longrightarrow
\mathbb{R}$, $k=\overline{1,4}$, of the expansion
\begin{equation}\label{mon-funk}
\Phi(\zeta)=U_{1}(x,y)\,e_1+U_{2}(x,y)\,ie_1+
U_{3}(x,y)\,e_2+U_{4}(x,y)\,ie_2
\end{equation}
satisfies also the equation (\ref{big-eq}), i.e. $U_{k}$ is a {\em
biharmonic function} in the domain $D$.

It is proved in the paper \cite{UMZ_PL_G_2009} that a monogenic
function $\Phi : D_{\zeta}\longrightarrow \mathbb{B}$ has
derivatives $\Phi^{(n)}(\zeta)$ of all orders in the domain
$D_{\zeta}$ and, consequently, satisfies the two-dimensional
biharmonic equation (\ref{big-eq}). In the papers
\cite{UMZ_PL_G_2009,Cont_2011} it was also proved such a fact that
every biharmomic function $U_{1}(x,y)$ in a bounded simply
connected domain $D$  is the first component of the expansion
(\ref{mon-funk}) of monogenic function $\Phi :
D_{\zeta}\longrightarrow \mathbb{B}$ determined in an explicit
form.

Basic analytic properties of monogenic functions in a biharmonic
plane are similar to properties of holomorphic functions of the
complex variable. More exactly, analogues of the Cauchy integral
theorem and integral formula, the Morera theorem, the uniqueness
theorem, the Taylor and Laurent expansions are established in the
paper \cite{Gr-Pl_Dop2009}.

\section{Statement of problem}
Consider the following boundary value problem: to find a monogenic
function $\Phi : D_{\zeta }\longrightarrow\mathbb{B}$ when values
of two components of the expansion (\ref{mon-funk}) are given on
the boundary $\partial D_{\zeta }$, i.e. the following boundary
conditions are satisfied:
$$U_{k}(x,y)=u_{k}(\zeta)\,,\quad
U_{m}(x,y)=u_{m}(\zeta)\qquad\forall\, \zeta \in
\partial D_{\zeta}\,,\quad 1\le k<m\le
4\,,$$ where $u_{k}$ and $u_{m}$ are given functions.

This problem was first considered by V.F. Kovalev \cite{Kov} and
was called as the {\it biharmonic Schwartz problem} because it is
analogous in a certain sense to the classical Schwartz problem on
finding an analytic function of complex variable when values of
its real part are given on the boundary of domain.

Let us analyse the main ideas of the paper \cite{Kov}. V.F.
Kovalev proved that all biharmonic Schwartz problem are reduced to
the main three problems:

1)  for $k=1$ and $m=2$ (we shall call it by the (1-2)-problem);

2)  for $k=1$ and $m=3$ (we shall call it by the (1-3)-problem);

3)  for $k=1$ and $m=4$ (we shall call it by the (1-4)-problem).

Some relations between biharmonic Schwartz problems and problems
of the theory of elasticity are described in the paper \cite{Kov}.
In particular, it is shown that the main biharmonic problem (see,
for example, \cite[p. 202]{Sm-3}) on finding a biharmonic function
$U : D \longrightarrow\mathbb{R}$ with given limiting values of
its partial derivatives $\partial U/
\partial x$ and $\partial U/ \partial y$ on the boundary $\partial D$ can be
reduced to the (1-3)-problem.

If certain natural conditions are satisfied, then the
(1-4)-problem is solved in \cite{Kov} in an explicit form.
Further, the (1-2)-problem and the (1-3)-problem are reduced to
integro-differential equations with using a conformal mapping of a
domain $D_{\zeta}$ onto the half-plane
$\Pi^{+}:=\{\zeta=xe_1+ye_2:y>0\}$ and the (1-4)-problem as an
auxiliary problem. Note that V.F. Kovalev \cite{Kov} stated only a
sketch of solving biharmonic Schwartz problems, and he did not
investigate conditions of solvability of these  problems.

Having an intention  to reduce the (1-3)-problem to integral
equations, we study boundary properties of certain integral
representations of monogenic functions. Considering the
(1-3)-problem for a half-plane and for a disk under natural
assumptions on the functions $u_1$ and $u_3$, we obtain solutions
in explicit forms. Moreover, the (1-3)-problem is solvable
unconditionally for a half-plane but it is solvable for a disk if
and only if a certain condition is satisfied.


\section{Biharmonic Schwartz integral for a half-plane}

Let a function $u : \mathbb{R}\longrightarrow \mathbb{R}$ be
continuous and there exists a finite limit
\begin{equation} \label{lim_infty}
u(\infty):=\lim\limits_{t\to \infty}u(t)\,.
\end{equation}
Under assumptions that the modulus of continuity
$$\omega_{\mathbb R}(u, \varepsilon)=\sup\limits_{\tau_1, \tau_2 \in \mathbb{R} : |\tau_1-\tau_2|\le \varepsilon }\left|u(\tau_1)-u(\tau_2)\right| $$
and the local centered (with respect to the infinitely remote
point) modulus of continuity
$$\omega_{\mathbb R, \infty}(u, \varepsilon)=\sup\limits_{\tau \in \mathbb{R} : |\tau|\ge 1/\varepsilon }\left|u(\tau)-u(\infty)\right|$$
of the function $u$ satisfy the Dini conditions
\begin{equation} \label{usl_Dini}
\int\limits_{0}^{1}\frac{\omega_{\mathbb R}(u,
\eta)}{\eta}\,d\,\eta<\infty,
\end{equation}
 \begin{equation} \label{usl_Dini+infty}
\int\limits_{0}^{1}\frac{\omega_{\mathbb R,  \infty}(u,
\eta)}{\eta}\,d\,\eta<\infty,
\end{equation}
consider an integral
$$ S_{\Pi^+}[u](\zeta):= \frac{1}{\pi
i}\int\limits_{-\infty}^{+\infty}\frac{u(t)(1+t\zeta)}{(t^2+1)}(t-\zeta)^{-1}\,dt\quad
\forall\, \zeta\in \Pi^{+} $$ that we shall call the {\it
biharmonic Schwartz integral} for the half-plane $\Pi^{+}$.

Here and in what follows, all integrals along the real axis are
understood in the sense of their Cauchy principal values, i.e.
$$\int\limits_{-\infty}^{+\infty}g(t,\cdot)\,dt:=\lim_{N\to
+\infty}\int\limits_{-N}^{N}g(t,\cdot)\,dt\,,$$
$$\int\limits_{-\infty}^{+\infty}\frac{g(t,\cdot)}{t-\xi}\,dt:=\lim_{N\to
+\infty}\,\,\lim_{\varepsilon\to 0+0}
\left(\int\limits_{-N}^{\xi-\varepsilon}+\int\limits_{\xi+\varepsilon}^{N}\right)
\frac{g(t,\cdot)}{t-\xi}\,dt\,,\quad \xi \in \mathbb{R}\,.$$

The function $S_{\Pi^+}[u](\zeta)$ is the principal extension (see
\cite[p. 165]{Hil-Fil}) into the half-plane $\Pi^{+}$ of the
complex Schwartz integral
$$ S[u](z):= \frac{1}{\pi
i}\int\limits_{-\infty}^{+\infty}\frac{u(t)(1+tz)}{(t^2+1)(t-z)}\,dt\,,$$
which determines a holomorphic function in the half-plane
$\{z=x+iy:y>0\}$ of the complex plane $\mathbb{C}$ with the given
boundary values $u(t)$ of real part on the real line $\mathbb{R}$.
Furthermore, the equality
\begin{equation}\label{bSw+kSw}
S_{\Pi^+}[u](\zeta)=S[u](z)-\frac{y}{2\pi}\,\rho\,
\int\limits_{-\infty}^{\infty}\frac{u(t)}{(t-z)^2}d\,t
\qquad\forall\, \zeta=xe_1+e_2y\in \Pi^{+}
\end{equation}
holds, where $z=x+iy$ as well as in what follows.


We use the euclidian norm $\|a\|:=\sqrt{|z_1|^2+|z_2|^2}$ in the
algebra $\mathbb{B}$, where\, $a=z_1e_1+z_2e_2$ and $z_1, z_2\in
\mathbb{C}$.

The following lemma presents sufficient conditions for the
existence of boundary values of the biharmonic Schwartz integral
on the extended real line ${\mathbb R} \cup \{\infty\}$.

\vskip 2mm \lem\label{lem1} {\it If a function $u :
\mathbb{R}\longrightarrow \mathbb{R}$ has the finite limit
{\emph{(\ref{lim_infty})}} and  the condition
{\emph{(\ref{usl_Dini})}} is satisfied, then the equality
\begin{equation} \label{gr_big_Sw}
\lim\limits_{\zeta\to \xi,\,
\zeta\in\Pi^+}S_{\Pi^+}[u](\zeta)=u(\xi)+\frac{1}{\pi
i}\int\limits_{-\infty}^{\infty}\frac{u(t)}{t^2+1}\frac{1+t\xi}{t-\xi}\,dt
\quad \forall\,\xi \in \mathbb{R}
\end{equation}
is fulfilled. If, in addition, the function $u$ satisfies the
condition {\emph{(\ref{usl_Dini+infty})}}, then
\begin{equation} \label{gr_big_Sw+infty}
\lim\limits_{\|\zeta\|\to \infty,\,
\zeta\in\Pi^+}S_{\Pi^+}[u](\zeta)=u(\infty)-\frac{1}{\pi
i}\int\limits_{-\infty}^{\infty}u(t)\frac{t}{t^2+1}\,dt.
\end{equation}}

{\bf Proof.} In order to prove the equality (\ref{gr_big_Sw}) we
use the expression (\ref{bSw+kSw}) of the biharmonic Schwartz
integral. The second summand in the right-hand part of equality
(\ref{bSw+kSw}) tends to zero with $\zeta\to\xi\in \mathbb{R}$.
This statement follows from the equalities
$$
y\int\limits_{-\infty}^{\infty}\frac{u(t)}{(t-z)^2}d\,t
=y\int\limits_{-\infty}^{\infty}\frac{u(t)-u(x)}{(t-z)^2}d\,t=
y\int\limits_{x-2|y|}^{x+2|y|}\frac{u(t)-u(x)}{(t-z)^2}\,dt+$$
\begin{equation}\label{uu}
+y\left(\int\limits_{-\infty}^{x-2|y|}+\int\limits_{x+2|y|}^{\infty}\right)\frac{u(t)-u(x)}{(t-z)^2}\,dt=:
I_1+I_2
\end{equation}
and the relations
\begin{equation}\label{uu1}
|I_1|\le
|y|\int\limits_{x-2|y|}^{x+2|y|}\frac{|u(t)-u(x)|}{y^2}\,dt\le
4\,\omega_{\mathbb{R}}(u,2|y|)\to 0,\quad z\to\xi\,,
\end{equation}
\begin{equation}\label{uu2}
|I_2|\le
|y|\left(\int\limits_{-\infty}^{x-2|y|}+\int\limits_{x+2|y|}^{\infty}\right)\frac{|u(t)-u(x)|}{|t-x|^2}\,dt
\le
2|y|\int\limits_{2|y|}^{\infty}\frac{\omega_{\mathbb{R}}(u,\eta)}{\eta^2}\,d\eta\to
0,\quad z\to\xi\,.
\end{equation}

By virtue of the condition (\ref{usl_Dini}), the Schwartz integral
$S[u](z)$ has limiting values on the real line (it follows, for
example, from an appropriate result of the paper \cite{Gerus_78}
for the Cauchy type integral), hence
$$
S[u](z)\to u(\xi)+\frac{1}{\pi
i}\int\limits_{-\infty}^{\infty}\frac{u(t)}{t^2+1}\frac{1+t\xi}{t-\xi}\,dt\,,
\quad z=x+iy\to \xi\in \mathbb{R},\,\,\, y>0\,.
$$

Thus, the equality (\ref{gr_big_Sw}) is proved.

In order to prove the equality (\ref{gr_big_Sw+infty}) with using
the change of variables $t=-1/t_1$, $z=-1/z_1$ (see., for example,
\cite[p. 36]{Gaxov}) we rewrite the relation (\ref{bSw+kSw}) in
the form
\begin{equation} \label{bSw+infty}
S_{\Pi^+}[u](\zeta)=S[v](z_1)-\frac{1}{2\pi}\,\frac{z_1}{\bar
z_1}\,\rho\,{\rm Im}\,z_1
\int\limits_{-\infty}^{\infty}\frac{v(t_1)}{(t_1-z_1)^2}d\,t_1\,,
\end{equation}
where $v(t_1):=u(-1/t)$ and $\bar z_1:={\rm Re}\,z_1-i\,{\rm
Im}\,z_1$. By virtue of Lemma 1 of the paper \cite{P-01-12}, the
function $v$ satisfies a condition of the form (\ref{usl_Dini}).
Therefore, the equality (\ref {gr_big_Sw+infty}) can be obtained
from the equality (\ref {bSw+infty}) by passing to the limit when
$z_1\to 0 $, \,  ${\rm Im}\,z_1>0$ by analogy with the proof of
the equality (\ref{gr_big_Sw}).  The lemma is proved.

\section{(1-3)-problem for a half-plane}

Let continuous functions $u_1 : \mathbb{R} \longrightarrow
\mathbb{R}$ and $u_3 : \mathbb{R} \longrightarrow \mathbb{R}$ have
finite limits of the 
form (\ref{lim_infty}).

Consider the following (1-3)-problem for the half-plane $\Pi^+$:
to find a continuous function $\Phi :
\overline{\Pi^+}\longrightarrow\mathbb{B}$ which is monogenic in
the domain $\Pi^+$ and has the limit
$$\lim\limits_{\|\zeta\|\to \infty,\,\zeta\in
\Pi^+}\Phi(\zeta)=:\Phi(\infty)\in \mathbb{B}\,,$$
when values of the components $U_1$ and $U_3$ of the expansion
(\ref{mon-funk}) are given on the real line $\mathbb{R}$, i.e. the
following boundary conditions are satisfied:
$$U_{1}(\xi,0)=u_{1}(\xi)\,,\quad
U_{3}(\xi,0)=u_{3}(\xi)\qquad\forall\, \xi\in \mathbb{R}\,.$$

It follows from Lemma \ref{lem1} that the function
 \begin{equation} \label{sol_1-ch}
\Phi(\zeta)=S_{\Pi^+}[u_1](\zeta)\,e_1+S_{\Pi^+}[u_3](\zeta)\,e_2
\end{equation}
is a solution of (1-3)-problem for the half-plane $\Pi^+$ if the
functions $u_1$ and $u_3$ satisfy the same conditions as the
function $u$ in Lemma \ref{lem1}.

To describe all solutions of (1-3)-problem for $\Pi^{+}$, first
consider the homogeneous (1-3)-problem.


\vskip 2mm \lem\label{homog_(1-3)_probl} {\it All solutions $\Phi
: \overline{\Pi^+}\longrightarrow\mathbb{B}$ of the homogeneous
{\emph{(1-3)}}-problem for the half-plane $\Pi^{+}$ with zero data
$u_1=u_3\equiv 0$ are expressed in the form
\begin{equation}\label{sol_homog}
\Phi(\zeta)= a_1\,ie_1+a_2\,ie_2\,,
\end{equation}
where $a_1, a_2$ are any real constants.}

{\bf Proof.} We use an expression of monogenic function $\Phi$ via
two holomorphic functions $F$ and $F_0$ (see, for example,
\cite{UMZ_PL_G_2009}) in the form
\begin{equation}\label{pr_mon}
\Phi(\zeta)=F(z)e_1-\left(\frac{iy}{2}F'(z)-F_{0}(z)\right)\rho
\qquad\forall\,\zeta\in \Pi^{+}.
\end{equation}

Inasmuch as the function $F$ is holomorphic in the half-plane
$\{z=x+iy: y>0\}$ and continuous on the closed half-plane
$\{z=x+iy: y\ge 0\}$ and have the limit
$$\lim\limits_{z\to \infty,\,\Re z>0}F(z)=:F(\infty)\,,$$
the following equalities hold:
$$F(z)=\frac{1}{2\pi i}
\int\limits_{-\infty}^{\infty}\frac{F(t)}{t-z}\,dt+\frac{F(\infty)}{2}\,,$$
\begin{equation}\label{F'}
F'(z)=\frac{1}{2\pi
i}\int\limits_{-\infty}^{\infty}\frac{F(t)}{(t-z)^2}\,dt\qquad
\forall\,z\in \mathbb{C} : \Re z>0\,.
\end{equation}
Therefore, we obtain the equality
$$ \lim\limits_{z\to\xi,\, y>0}\,yF'(z)=0,\qquad \forall\,\xi\in\mathbb{R}$$
as the result of relations (\ref{uu}) --- (\ref{uu2}).

In the same way with using the change of variables $t=-1/t_1$,
$z=-1/z_1$ in the integral (\ref{F'}) we obtain the equality
$$ \lim\limits_{z\to\infty,\,y>0}\,yF'(z)=0.$$

Thus, taking into account the relation (\ref{rho}), we have
$$  \lim\limits_{\zeta\to\xi,\,\zeta\in \Pi^{+}}\,\Phi(\zeta)=\Bigl(F(\xi)+2F_0(\xi)\Bigr)e_1+2iF_0(\xi)e_2
\qquad\forall\,\xi\in\mathbb{R},$$
 and the homogeneous (1-3)-problem for $\Pi^{+}$ is reduced to finding holomorphic
functions $F$, $F_0$ by solving two classical Schwartz
problem for the half-plane $\{z=x+iy: y>0\}$ with the following
boundary conditions:
$$\Re \Bigl(F(\xi)+2F_0(\xi)\Bigr)=0, \qquad \Re
\Bigl(2iF_0(\xi)\Bigr)=0\qquad\forall\,\xi\in\mathbb{R}.$$
 In such a way we obtain $F_{0}(z)\equiv a_2/2$ and $F(z)\equiv -a_2+ia_1$,
where $a_1, a_2$ are any real constants.

Finally, substituting the obtained functions $F$, $F_0$ into the
equality (\ref{pr_mon}) and taking into account the relation
(\ref{rho}), we obtain the equality (\ref{sol_homog}). The lemma
is proved.

In the following theorem we establish the formula of solutions of
the (1-3)-problem for the half-plane $\Pi^{+}$.

 \vskip 2mm \theor \label{Z_1-3}
 {\em Let the functions $u_1 : \mathbb{R} \longrightarrow
\mathbb{R}$ and $u_3 : \mathbb{R} \longrightarrow \mathbb{R}$ have
finite limits of the form {\emph{(\ref{lim_infty})}} and satisfy
conditions of the form {\emph{(\ref{usl_Dini})}} and
{\emph{(\ref{usl_Dini+infty})}}. Then the general solution of
{\emph{(1-3)}}-problem for $\Pi^{+}$ is expressed in the form
 \begin{equation} \label{sol_1-3}
 \Phi(\zeta)=S_{\Pi^+}[u_1](\zeta)\,e_1+S_{\Pi^+}[u_3](\zeta)\,e_2 +a_1\,ie_1+a_2\,ie_2,
 \end{equation}
where $a_1, a_2$ are any real constants.}

{\bf Proof.} It is obvious that the formula (\ref{sol_1-3})
represents the solution of the (1-3)-problem for $\Pi^{+}$ as the
sum of particular solution (\ref{sol_1-ch}) and the general
solution (\ref{sol_homog}) of the homogeneous (1-3)-problem. The
theorem is proved.

\section{A biharmonic analogue of Schwartz integral for a disk}

In what follows, $D_{\zeta}:=\{\zeta=xe_1+ye_2 : \|\zeta\|\le 1\}$
is the unit disk in the biharmonic plane $\mu$ and $D:=\{z=x+iy :
|z|\le 1\}$ is the unit disk in the complex plane $\mathbb{C}$.

For a continuous function $u : \partial D_{\zeta} \longrightarrow
\mathbb{R}$, by $\widehat{u}$ we denote the function defined on
the unit circle $\partial D$ of the complex plane $\mathbb{C}$ by
the equality $\widehat{u}(z)=u(\zeta)$ for all $z\in \partial D$.

Consider the integral
\begin{equation}\label{big_int_Sw_c}
S_{D_{\zeta}}[u](\zeta):= \frac{1}{2\pi i}\int\limits_{\partial
D_{\zeta}}
u(\tau)(\tau+\zeta)(\tau-\zeta)^{-1}\,\tau^{-1}\,d\tau\qquad
\forall\, \zeta\in D_{\zeta}
\end{equation}
that is an analogue of the complex Schwartz integral
\begin{equation}\label{int_Sw_c}
\frac{1}{2\pi i}\int\limits_{\partial
D}\widehat{u}(t)\,\frac{t+z}{t-z}\,\frac{dt}{t} \qquad \forall\,
z\in D
\end{equation}
which determines a holomorphic function in the disk $D$ with the
given boundary values $\widehat{u}(t)$ of real part on the circle
$\partial D$.

Consider also singular integrals which are understood in the sense
of their Cauchy principal values, i.e.
$$ S_{\partial D_{\zeta}}[u](\zeta):=\frac{1}{2\pi i}\,\lim_{\varepsilon\to 0+0} \int\limits_{\{\tau\in\partial
D_{\zeta}\, :\, \|\tau-\zeta\|\ge\varepsilon\} }
u(\tau)(\tau+\zeta)(\tau-\zeta)^{-1}\,\tau^{-1}\,d\tau\qquad
\forall\, \zeta\in\partial D_{\zeta}\,,$$
$$S_0[\widehat{u}](z):=\frac{1}{2\pi i}\,\lim_{\varepsilon\to 0+0} \int\limits_{\{t\in\partial
D\, : \, |t-z|\ge\varepsilon\}
}\widehat{u}(t)\,\frac{t+z}{t-z}\,\frac{dt}{t} \qquad \forall\,
z\in\partial D\,.$$

It can be proved in a similar way as in the complex plane (see,
for example, \cite{Gerus_78,Gaxov}) that if the modulus of
continuity
$$\omega(u,
\varepsilon)=\sup\limits_{\tau_1, \tau_2 \in\partial D_{\zeta} :
\|\tau_1-\tau_2\|\le \varepsilon} \|u(\tau_1)-u(\tau_2)\|
$$ of the function $u : \partial D_{\zeta} \longrightarrow
\mathbb{R}$ satisfies the Dini condition
\begin{equation} \label{usl_Dini_c}
\int\limits_{0}^{1}\frac{\omega(u, \eta)}{\eta}\,d\,\eta<\infty\,,
\end{equation}
then the integral (\ref{big_int_Sw_c}) has limiting values on
$\partial D_{\zeta}$ which are expressed by the formula
\begin{equation} \label{gr_big_Sw_c}
\lim\limits_{\xi\to\zeta,\, \xi\in D_{\zeta}}
S_{D_{\zeta}}[u](\xi)=u(\zeta)\,e_1+S_{\partial
D_{\zeta}}[u](\zeta) \qquad \forall\, \zeta\in\partial
D_{\zeta}\,.
\end{equation}


\vskip 2mm \lem\label{expr_big_Sw_c} {\it If the function $u :
\partial D_{\zeta} \longrightarrow \mathbb{R}$ satisfies the
condition (\ref{usl_Dini_c}), then the following equality holds:
$$
S_{\partial D_{\zeta}}[u](\zeta)=
S_0[\widehat{u}](z)\,e_1-\left(\frac{y}{2
\pi}\int\limits_{\partial D}\frac{\widehat{u}(t)}{t^2}\,dt \right)
\,(e_1+ie_2)+$$
\begin{equation}\label{big_Sw_c+int_Sw_c}
+\left(\frac{x}{2\pi}\int\limits_{\partial D}
\frac{\widehat{u}(t)}{t^2}\,dt+\frac{1}{2\pi}\int\limits_{\partial
D} \frac{\widehat{u}(t)}{t^3}\,dt\right)(e_2-ie_1)
\quad\forall\,\zeta\in \partial D_{\zeta}.
\end{equation} }

{\bf Proof.} Let $\zeta\in\partial D_{\zeta}$,
$\tau:=t_1e_1+t_2e_2$ and $t:=t_1+it_2$, where
$t_1,t_2\in\mathbb{R}$. We have the equalities
$\zeta=ze_1-\frac{iy}{2}\rho$, $\tau=te_1-\frac{it_2}{2}\rho$ and
$d\tau=e_1dt-\frac{i}{2}\rho\,dt_2$ for $\tau\in\partial
D_{\zeta}$. Taking also into account the equalities (see
\cite{UMZ_PL_G_2009})
$$(\tau-\zeta)^{-1}=\frac{1}{t-z}\,e_1+\frac{i(t_2-y)}{(t-z)^{2}}\,\rho\,,\qquad
\tau^{-1}=\frac{1}{t}\,e_1+\frac{it_2}{t^2}\,\rho\,,$$
 we obtain the following equality:
$$S_{\partial D_{\zeta}}[u](\zeta)=
S_0[\widehat{u}](z)\,e_1+\frac{1}{2\pi i}\int\limits_{\partial D}
\frac{\widehat{u}(t)(xt_2-t_1y)}{t^2(t-z)}\,dt\left(\frac{i\rho}{2}\right)+$$
$$+\frac{1}{2\pi i}\int\limits_{\partial D}
\frac{\widehat{u}(t)(t+z)((t_2-y)\,dt_1-(t_1-x)\,dt_2)}{t(t-z)^{2}}\left(\frac{i\rho}{2}\right)=:$$

\begin{equation}\label{felicita'}
=:S_0[\widehat{u}](z)\,e_1+\frac{i\rho}{2}\,I_3+
\frac{i\rho}{2}\,I_4.
\end{equation}

Transforming the integrals $I_3$ and $I_4$ with using the change
of variables $t=\exp(i\theta)$, $z=\exp(i\theta_z)$, we obtain the
equalities
$$I_3=I_4=\frac{-x+iy}{4 \pi}\int\limits_{\partial D} \frac{\widehat{u}(t)}{t^2}\,dt-\frac{1}{4 \pi}\int\limits_{\partial D} \frac{\widehat{u}(t)}{t^3}\,dt\,.$$

Substituting obtained expressions for $I_3$ and $I_4$ into the
equality (\ref{felicita'}) and taking into account the equality
$i\rho/2=ie_1-e_2$  which follows from the equality (\ref{rho}),
we obtain the equality (\ref{big_Sw_c+int_Sw_c}). The lemma is
proved.

\section{(1-3)-problem for a disk}

Let the functions $u_1 : \partial D_{\zeta} \longrightarrow
\mathbb{R}$ and $u_3 : \partial D_{\zeta} \longrightarrow
\mathbb{R}$ satisfy the condition of the form (\ref{usl_Dini_c}).

Consider the following (1-3)-problem for the unit disk
$D_{\zeta}$: to find a continuous function $\Phi :
\overline{D_{\zeta}} \longrightarrow\mathbb{B}$ which is monogenic
in the disk $D_{\zeta}$ when values of the components $U_1$ and
$U_3$ of the expansion (\ref{mon-funk}) are given on the circle
$\partial D_{\zeta}$, i.e. the following boundary conditions are
satisfied:
\begin{equation} \label{u_k,m_c}
U_{1}(x,y)=u_{1}(\zeta)\,,\quad
U_{3}(x,y)=u_{3}(\zeta)\qquad\forall\, \zeta \in \partial
D_{\zeta}.
\end{equation}

It follows from a biharmonic analogue of the Cauchy integral
formula (see Theorem 4 from the paper \cite{Gr-Pl_Dop2009}) that
the solution $\Phi(\zeta)$ of 
(1-3)-problem for $D_{\zeta}$
can be represented as a biharmonic Cauchy integral and,
consequently, it can also be represented in the form
\begin{equation} \label{form-sol}
\Phi(\zeta)= \frac{1}{2\pi i}\int\limits_{\partial D_{\zeta}}
\varphi(\tau)(\tau+\zeta)(\tau-\zeta)^{-1}\,\tau^{-1}\,d\tau\qquad
\forall\, \zeta\in D_{\zeta}\,,
 \end{equation}
 where $\varphi : \partial D_{\zeta}\longrightarrow
\mathbb{B}$ is a certain continuous function.

We shall find solutions of 
(1-3)-problem for $D_{\zeta}$ in the class $\mathcal{M}$ of
functions represented in the form (\ref{form-sol}), where the
function $\varphi : \partial D_{\zeta}\longrightarrow \mathbb{B}$
satisfies the condition of the form (\ref{usl_Dini_c}).  The
solvability of (1-3)-problem for $D_{\zeta}$ in the class
$\mathcal{M}$ is described in the following theorem.

 \vskip 2mm \theor \label{Z_1-3_d}
 {\em Let the functions $u_1 : \partial D_{\zeta} \longrightarrow
\mathbb{R}$ and $u_3 : \partial D_{\zeta} \longrightarrow
\mathbb{R}$ satisfy the condition of the form (\ref{usl_Dini_c}).
Then {\emph{(1-3)}}-problem for $D_{\zeta}$ is solvable in the
class $\mathcal{M}$ if and only if the following condition is
satisfied:
\begin{equation}\label{dop_usl_d}
\int\limits_{\partial D_{\zeta}} u_1(\zeta)\,dx+u_3(\zeta)\,dy =0.
\end{equation}
 The general solution is expressed in the form
 \begin{equation} \label{sol_1-3_d} \Phi(\zeta)=
S_{D_{\zeta}}[u_1](\zeta)\,e_1+S_{D_{\zeta}}[u_3](\zeta)\,e_2+b\zeta+b_1e_1+b_2e_2+
i(a\zeta +a_1\,e_1+a_2\,e_2),
 \end{equation}
where
$$b:=-\frac{1}{2\pi}\Re\int\limits_{\partial D}\frac{\widehat{u}_{3}(t)}{t^2}\,dt-\frac{1}{2\pi}\Im\int\limits_{\partial D}\frac{\widehat{u}_{1}(t)}{t^2}\,dt\,,$$
$$b_1:=-\frac{1}{2\pi}\Re\int\limits_{\partial D}\frac{\widehat{u}_{3}(t)}{t^3}\,dt-\frac{1}{2\pi}\Im\int\limits_{\partial D}\frac{\widehat{u}_{1}(t)}{t^3}\,dt\,,$$
$$b_2:=\frac{1}{2\pi}\Im\int\limits_{\partial D}\frac{\widehat{u}_{3}(t)}{t^3}\,dt-\frac{1}{2\pi}\Re\int\limits_{\partial D}\frac{\widehat{u}_{1}(t)}{t^3}\,dt\,,$$
and $a, a_1, a_2$ are any real constants.}

{\bf Proof.} Let $\Phi\in\mathcal{M}$ be a solution of 
(1-3)-problem for $D_{\zeta}$. Taking into account the equalities
(\ref{gr_big_Sw_c}) and (\ref{big_Sw_c+int_Sw_c}), we make a
conclusion that the limiting values of the function $\Phi$ on the
boundary $\partial D_{\zeta}$ can be represented in the form
\begin{equation} \label{b-v-sol_1-3_d}
\Phi(\zeta)= F_1(z)\,e_1+ F_2(z)\,e_2 +c_1x+c_2y+c_3 \qquad
\forall\, \zeta\in
\partial D_{\zeta}\,,
 \end{equation}
 where $F_1 : \overline{D}\longrightarrow\mathbb{C}$,
$F_2 : \overline{D}\longrightarrow \mathbb{C}$ are certain
continuous functions holomorphic in the disk $D$, and $c_1, c_2,
c_3\in\mathbb{B}$. Moreover, the functions $F_1$, $F_2$ can be
represented by complex Schwartz integrals of the form
(\ref{int_Sw_c}). Therefore, the equality (\ref{b-v-sol_1-3_d})
can be rewritten as
\begin{equation} \label{b-v-sol_1-3_d-2}
\Phi(\zeta)= g_1(\zeta)\,e_1+S_{\partial D_{\zeta}}[g_1](\zeta)
+e_2\Bigl(g_3(\zeta)\,e_1+S_{\partial D_{\zeta}}[g_3](\zeta)\Bigr)
+h_1x+h_2y+h_3 \qquad \forall\, \zeta\in
\partial D_{\zeta}\,,
 \end{equation}
where $g_1 : \partial D_{\zeta} \longrightarrow \mathbb{R}$, $g_3
: \partial D_{\zeta} \longrightarrow \mathbb{R}$ are certain
functions satisfying the condition of the form (\ref{usl_Dini_c}),
and $h_1, h_2, h_3\in\mathbb{B}$.

It follows from a biharmonic analogue of the Cauchy integral
theorem (see Theorem 4 from the paper \cite{Gr-Pl_Dop2009}) that
for the function (\ref{b-v-sol_1-3_d-2}) the following equalities
must be fulfilled:
$$\int\limits_{\partial D_{\zeta}}\Phi(\zeta)\,d\zeta=
\int\limits_{\partial D_{\zeta}} (h_1x+h_2y)\,d\zeta=0\,.$$ The
last equality is feasible if and only if\,\, $h_2=h_1e_2$.

Denote $U_{k}[\Phi(\zeta)]:=U_{k}(x,y)$, where $U_{k}(x,y)$,
$k=\overline{1,4}$, are components of the expansion
(\ref{mon-funk}).

Our strategy is to find the functions $g_1$ and $g_3$ with using
the method of indefinite coefficients and the boundary conditions
(\ref{u_k,m_c}) rewritten in the form of two equations with
respect to sought-for functions $g_1$, $g_3$ and unknown
hypercomplex numbers $h_1, h_3$:
\begin{equation}\label{1_3_sys}
g_{k}(\zeta)+U_{k}\left[S_{\partial
D_{\zeta}}[g_1](\zeta)\,e_1+S_{\partial
D_{\zeta}}[g_3](\zeta)\,e_2\right]+xU_{k}[h_1]+yU_{k}[h_1e_2]+U_{k}[h_3]=u_{k}(\zeta)\qquad
\forall \zeta\in \partial D_{\zeta},\, k=1,3.
\end{equation}

Using the equality (\ref{big_Sw_c+int_Sw_c}), we obtain the
following equalities:
$$U_{1}\left[S_{\partial D_{\zeta}}[g_1](\zeta)\,e_1+
S_{\partial D_{\zeta}}[g_3](\zeta)\,e_2\right]=\frac{1}{2\pi}\left( \Im \int\limits_{\partial D}\frac{\widehat{g}_1(t)}{t^2}\,dt+
\Re\int\limits_{\partial D}
\frac{\widehat{g}_3(t)}{t^2}\,dt\right)x+$$
\begin{equation} \label{u_1}
+\frac{1}{2\pi}\left(-\Re\int\limits_{\partial D}
\frac{\widehat{g}_1(t)}{t^2}\,dt+ \Im\int\limits_{\partial D}
\frac{\widehat{g}_1(t)}{t^3}\,dt\right)y+
\frac{1}{2\pi}\Im\int\limits_{\partial D}
\frac{\widehat{g}_1(t)}{t^2}\,dt+
\frac{1}{2\pi}\Re\int\limits_{\partial D}
\frac{\widehat{g}_3(t)}{t^3}\,dt,
\end{equation}
$$U_{3}\left[S_{\partial D_{\zeta}}[g_1](\zeta)\,e_1+S_{\partial D_{\zeta}}[g_3](\zeta)\,e_2\right]=
\frac{1}{2\pi}\left(\Re\int\limits_{\partial D}\frac{\widehat{g}_1(t)}{t^2}\,dt
-\Im\int\limits_{\partial D}
\frac{\widehat{g}_3(t)}{t^2}\,dt\right)x+$$
\begin{equation}\label{u_3}
+\frac{1}{2\pi}\left(\Im\int\limits_{\partial D}
\frac{\widehat{g}_1(t)}{t^2}\,dt+ \Re\int\limits_{\partial D}
\frac{\widehat{g}_3(t)}{t^2}\,dt\right)y
+\frac{1}{2\pi}\Re\int\limits_{\partial D}
\frac{\widehat{g}_1(t)}{t^3}\,dt
-\frac{1}{2\pi}\Im\int\limits_{\partial D}
\frac{\widehat{g}_3(t)}{t^3}\,dt.
\end{equation}

Now, it follows from the equalities (\ref{1_3_sys}) ---
(\ref{u_3}) that the functions $g_1$, $g_3$ can be expressed in
the form
\begin{equation}\label{g_k}
g_{k}(\zeta)=u_{k}(\zeta)+a_{k,1}x+a_{k,2}y+a_{k,0}\qquad \forall
\zeta\in \partial D_{\zeta},
\end{equation}
where unknown coefficients $a_{k,m}$ are real numbers for\,\,
$k=1, 3$\,\, and\,\, $m=0, 1, 2$\,.

Denote
$$\begin{array}{l}
\displaystyle
 A_{k}:=\frac{1}{2\pi}\Re\int\limits_{\partial D}\frac{\widehat{u}_{k}(t)}{t^2}\,dt,\quad
B_{k}:=\frac{1}{2\pi}\Im\int\limits_{\partial D}\frac{\widehat{u}_{k}(t)}{t^2}\,dt,\\[6mm]
\displaystyle
 C_{k}:=\frac{1}{2\pi}\Re\int\limits_{\partial D}\frac{\widehat{u}_{k}(t)}{t^3}\,dt,\quad
D_{k}:=\frac{1}{2\pi}\Im\int\limits_{\partial
D}\frac{\widehat{u}_{k}(t)}{t^3}\,dt, \qquad k=1,3.
\end{array}
$$

Substituting the expressions (\ref{g_k}) into the equations
(\ref{1_3_sys}) and taking into account the equalities
$$\frac{1}{2\pi}\int\limits_{\partial D}\frac{x}{z^2}\,dz=\frac{i}{2}\,,\qquad
\frac{1}{2\pi}\int\limits_{\partial
D}\frac{y}{z^2}\,dz=\frac{1}{2}\,,\qquad
 \int\limits_{\partial D}\frac{x}{z^3}\,dz=\int\limits_{\partial
 D}\frac{y}{z^3}\,dz=0\,,$$
  $U_1[h_1\,e_2]=U_{3}[h_1],\quad U_{3}[h_1\,e_2]=U_{1}[h_1]-2U_{4}[h_1]$,
we obtain the relations
$$
\begin{array}{l}
\displaystyle\left(\frac{3}{2}\,a_{1,1}+\frac{1}{2}\,a_{3,2}+A_3+B_1+U_1[h_1]\right)x+
\left(\frac{1}{2}\,a_{1,2}+\frac{1}{2}\,a_{3,1}-A_1+B_3+U_3[h_1]\right)y+\\[4mm]
\hspace*{70mm}+a_{1,0}+C_3+D_1+U_1[h_3]=0,\\[4mm]
\displaystyle\left(\frac{1}{2}\,a_{1,2}+\frac{1}{2}\,a_{3,1}+A_1-B_3+U_3[h_1]\right)x+
\left(\frac{1}{2}\,a_{1,1}+\frac{3}{2}\,a_{3,2}+A_3+B_1+U_1[h_1]-2U_4[h_1]\right)y+\\[4mm]
\hspace*{70mm}+a_{3,0}+C_1-D_3+U_3[h_3]=0.
\end{array}
$$

Consequently, we have a system of six equations with six real
unknowns $a_{k,m}$ (with $k=1, 3$ and\break $m=0, 1, 2$) and two
hypercomplex parameters $h_1, h_3$:
\begin{equation} \label{syst1}
\begin{array}{l}
a_{1,0}=-C_3-D_1-U_1[h_3],\\[2mm]
a_{3,0}=-C_1+D_3-U_3[h_3],\\[2mm]
a_{1,2}+a_{3,1}=2A_1-2B_3-2U_3[h_1],\\[2mm]
a_{1,2}+a_{3,1}=-2A_1+2B_3-2U_3[h_1],\\[2mm]
3a_{1,1}+a_{3,2}=-2A_3-2B_1-2U_1[h_1],\\[2mm]
a_{1,1}+3a_{3,2}=-2A_3-2B_1-2U_1[h_1]+4U_4[h_1].
\end{array}
\end{equation}

It is obvious that this system is solvable if and only if
$A_1-B_3=B_3-A_1$, i.e. $A_1=B_3$ that is equivalent to the
condition (\ref{dop_usl_d}). If this condition is satisfied, then
the general solution of the system (\ref{syst1}) contains an
arbitrary real number $a_{1,2}$ and is of the form:
\begin{equation} \label{sol-syst1}
\begin{array}{l}
a_{1,0}=-C_3-D_1-U_1[h_3],\\[2mm]
a_{3,0}=-C_1+D_3-U_3[h_3],\\[2mm]
a_{3,1}=-a_{1,2}-2U_3[h_1],\\[2mm]
\displaystyle a_{1,1}=-\frac{1}{2}\,A_3-\frac{1}{2}\,B_1-\frac{1}{2}\,U_1[h_1]-\frac{1}{2}\,U_4[h_1],\\[3mm]
\displaystyle
a_{3,2}=-\frac{1}{2}\,A_3-\frac{1}{2}\,B_1-\frac{1}{2}\,U_1[h_1]+\frac{3}{2}\,U_4[h_1].
\end{array}
\end{equation}

Inasmuch as\, $h_2=h_1e_2$, the function
\begin{equation} \label{sol_1-3_disk}
\Phi(\zeta)= S_{D_{\zeta}}[g_1](\zeta)\,e_1
+S_{D_{\zeta}}[g_3](\zeta)\,e_2+h_1\zeta+h_3 \qquad \forall\,
\zeta\in D_{\zeta}
 \end{equation}
 have the boundary values (\ref{b-v-sol_1-3_d-2}) and is the general solution of (1-3)-problem for $D_{\zeta}$ in the
class $\mathcal{M}$. In the formula (\ref{sol_1-3_disk}) the
functions $g_1, g_3$ are of the form (\ref{g_k}), where the
coefficients $a_{k,m}$ with\,\, $k=1, 3$\,\, and\,\, $m=0, 1, 2$\,
are determined by the equalities (\ref{sol-syst1}).

Finally, with using the equalities
$$S_{D_{\zeta}}[1](\zeta)=e_1,\quad S_{D_{\zeta}}[x](\zeta)=\frac{1}{2}\,(3e_1+ie_2)\zeta,
\quad S_{D_{\zeta}}[y](\zeta)=\frac{1}{2}\,(-3ie_1+e_2)\zeta
\qquad \forall\, \zeta\in D_{\zeta}$$ the formula
(\ref{sol_1-3_disk}) is reduced to the form (\ref{sol_1-3_d}). The
theorem is proved.\vskip 2mm

Let $E:=\{(x,y)\in\mathbb{R}^{2} : x^2+y^2<1\}$ be the unit disk
in the plane $\mathbb{R}^{2}$.

Consider the main biharmonic problem (see, for example, \cite[p.
202]{Sm-3}) on finding a biharmonic function $V : E
\longrightarrow\mathbb{R}$ with given limiting values of its
partial derivatives on the boundary $\partial E$:
\begin{equation}\label{osn_big_pr}
\begin{array}{l}
\displaystyle \lim\limits_{(x,y)\to(x_0,y_0),\, (x,y)\in E}
\frac{\partial
V(x,y)}{\partial x}=u_{1}(x_0,y_0),\\[4mm]
\displaystyle \lim\limits_{(x,y)\to(x_0,y_0),\, (x,y)\in E}
\frac{\partial V(x,y)}{\partial y}=u_{3}(x_0,y_0)\qquad \forall\,
(x_0,y_0) \in\partial E.
\end{array}
\end{equation}

It is shown in the paper \cite{Kov} that this problem can be
reduced to the (1-3)-problem. Indeed, consider in the domain
$D_{\zeta}$ a monogenic function
 $$\Phi_{1}(\zeta):=V(x,y)\,e_1+V_{2}(x,y)\,ie_1+V_{3}(x,y)\,e_2+V_{4}(x,y)\,ie_2$$
for which $U_{1}[\Phi_1(\zeta)]=V(x,y)$. It follows from the
Cauchy--Riemann condition (\ref{usl_K_R}) with $\Phi=\Phi_1$ that
$\partial V_{3}(x,y)/\partial x=\partial V(x,y)/\partial y$.
Therefore,
$$\Phi_{1}'(\zeta)=\frac{\partial V (x,y)}{\partial x}\,e_1+ \frac{\partial V_{2}(x,y)}{\partial x}\,ie_1+
\frac{\partial V(x,y)}{\partial y}\,e_2+ \frac{\partial
V_{4}(x,y)}{\partial x}\,ie_2\,,$$ and the main biharmonic problem
with the boundary conditions (\ref{osn_big_pr}) can be reduced to
(1-3)-problem on finding a monogenic function
$\Phi(\zeta):=\Phi_{1}'(\zeta)$ in the domain $D_{\zeta}$.

Let us note that in this case the equality (\ref{dop_usl_d}) is
the well known necessary and sufficient condition for the
solvability of the main biharmonic problem with the boundary
conditions (\ref{osn_big_pr}), see, for example, \cite[p.
202]{Sm-3}.

Let us also note that the solvability of the (1-3)-problem is
analogous to the solvability of Schwartz-type boundary value
problems for polyanalytic functions and for analytic functions of
several complex variables. In the paper \cite{Begehr_bi-pol} the
Schwartz-type boundary value problem for the polyanalytic equation
in a half-plane is solved in an explicit form without any
complementary conditions. In the paper \cite{Begehr_disk}
Schwartz-type boundary value problems for analytic functions of
several complex variables and for the inhomogeneous
Cauchy--Riemann system in a polydisc are considered, ibid.
solvability conditions for appropriate boundary value problems are
obtained.

\section{Acknowledgements}

The authors are grateful  to Professor Massimo Lanza de
Cristoforis for the invitation and given possibilities for
research during the {\em{Minicorsi di Analisi Matematica}}~--~2011
at the University of Padova.


\vspace{10mm}

Authors:\\[2mm]
{\bf S. V. Gryshchuk}, {University of Padova, Department of Pure
and Applied Mathematics, Italy, 35121, Padova, Via Trieste 63,\\
Phone (office): (0039) 049 827 1460,  E-mail:
serhii@math.unipd.it, gr\_sergey\_v@mail.ru}\\[2mm]
{ \bf S. A. Plaksa}, {Institute of Mathematics of the National
Academy of Sciences of Ukraine, Ukraine, 01601, Kiev,
Tereshchenkivska Str. 3,\\ Phone (office): (38044) 234 51 50,
E-mail: plaksa@imath.kiev.ua}

\end{document}